\documentclass[11pt]{amsart}

\usepackage{amsmath}
\usepackage{fullpage}
\usepackage{xspace}
\usepackage[psamsfonts]{amssymb}
\usepackage[latin1]{inputenc}
\usepackage{graphicx,color}
\usepackage[curve]{xypic} 
\usepackage{hyperref}
\usepackage{graphicx}
\usepackage[greek,english]{babel}

\usepackage{amsmath}%
\usepackage{amsthm}%
\usepackage{amscd}
\usepackage{amsfonts}%
\usepackage{amssymb}%
\usepackage{graphicx}

\usepackage{mathrsfs}

\usepackage{tikz}
\usetikzlibrary{matrix,arrows}

\usepackage{tikz-cd}

\newtheorem{theorem}{Theorem}[section]

\newtheorem{lemma}[theorem]{Lemma}

\newtheorem{proposition}[theorem]{Proposition}

\theoremstyle{remark}

\numberwithin{equation}{section}

\newcommand{\pfrak}{\mathfrak{p}}
\newcommand{\Pfrak}{\mathfrak{P}}

\newcommand{\Dcal}{\mathscr{D}}

\newcommand{\Fcal}{\mathscr{F}}

\newcommand{\Pcal}{\mathscr{P}}

\newcommand{\Z}{\mathbb{Z}}
\newcommand{\C}{\mathbb{C}}

\newcommand{\F}{\mathbb{F}}
\newcommand{\Q}{\mathbb{Q}}

\newcommand{\N}{\mathbb{N}}

\newcommand{\Gal}{\mathrm{Gal}}

\makeatletter
\@namedef{subjclassname@2020}{%
  \textup{2020} Mathematics Subject Classification}
\makeatother

  \DeclareFontFamily{U}{wncy}{}
    \DeclareFontShape{U}{wncy}{m}{n}{<->wncyr10}{}
    \DeclareSymbolFont{mcy}{U}{wncy}{m}{n}
    \DeclareMathSymbol{\Sha}{\mathord}{mcy}{"58}

\begin{document}
\title[]{Diophantine sets over global function fields}

\author{Hector Pasten}
\address{ Departamento de Matem\'aticas,
Pontificia Universidad Cat\'olica de Chile, and CENIA. Chile}
\email[H. Pasten]{hector.pasten@uc.cl}%

\author{Xavier Vidaux}
\address{ Departamento de Matem\'aticas,
Universidad de Concepci\'on. Chile}
\email[X. Vidaux]{xvidaux@udec.cl}%


\thanks{H.P. was supported by ANID Fondecyt Regular grant 1230507 from Chile. X.V. was supported by Proyecto 2025001477INV from Universidad de Concepci\'on.}
\date{\today}

\begin{abstract} Let $q\equiv 3\bmod 4$ be a prime power and let $K$ be a finite separable extension of the field of rational functions $\mathbb F_q(t)$, with constant field $\mathbb F_q$. We show that $\mathbb F_q[t]$ is Diophantine in $\mathbb F_q(t)$ and that the latter is Diophantine in $K$. We deduce from work of Demeyer that $K$ has the DPRM property: Diophantine sets over $K$ are the same as listable sets.
\end{abstract}

\maketitle



\section{Introduction}

Hilbert's tenth problem asked for an algorithm able to decide solvability in $\Z$ of Diophantine equations. The negative solution given by Davis, Putnam, Robinson, and Matiyasevich \cite{DPR, Matiyasevich} actually gives more than undecidability: Diophantine sets in $\Z^n$ are precisely the listable sets. 

More generally, we say that a listable ring $R$ has the DPRM property if the Diophantine sets (over $R$) in cartesian powers of $R$ are the same as the listable sets. See \cite{PastenDPRM}. So, $\Z$ has the DPRM property.

We are concerned with the case of global function fields. Let $q$ be a prime power, $k=\F_q$, and $F=k(t)$ with $t$ an independent variable. By \emph{global function field} we mean a finite extension $K/F$ which is separable and with constant field $k$.

First, Pheidas \cite{PheidasH10ratfn} gave a negative answer to Hilbert's tenth problem on $k(t)$ where $k$ is a finite field of odd characteristic. Characteristic $2$ was done by Videla \cite{Videla} and then the method was extended to cover general global function fields \cite{Shlapentokh}. This method, however, only yields undecidability and it does not provide a classification of Diophantine sets. We show the following results:
\begin{theorem}\label{ThmA}
Assume $q\equiv 3\bmod 4$. Then $k[t]$ is Diophantine in $F=k(t)$.
\end{theorem}
\begin{theorem}\label{ThmB}
Assume $q\equiv 3\bmod 4$ and let $K/F$ be a global function field. Then $F$ is Diophantine in $K$.
\end{theorem}

From these two results, by using work of Demeyer \cite{Demeyer} one immediately obtains:

\begin{theorem}\label{ThmC}
Assume $q\equiv 3\bmod 4$ and let $K/F$ be a global function field. Then $K$ has the DPRM property.
\end{theorem}

At the core of our proofs is the notion of a \emph{Pheidas predicate}, modeled after Pheidas's work \cite{Pheidas}. This notion has been central in all recent developments around Hilbert's tenth problem over $\Q$, see \cite{OpenAI, Glazer, Burungale, BhargavaBurungale, KPV}. After \cite{OpenAI}, the standard notation for a Pheidas predicate became $\Phi$.

Let us make some remarks.

\begin{itemize}

\item Our construction of the Pheidas predicate follows what is now standard after \cite{OpenAI}: points on two elliptic curves with a suitable squareclass condition. It is here where $q\equiv 3\bmod 4$ is used.

\item Our results sharply contrast those of Koll\'ar over $\C$, see \cite{Kollar}: in our case polynomials are Diophantine in rational functions, and for the function field extension $K/F$ we have that $F$ is Diophantine in $K$. 

\item Now we know that $\Z$ is Diophantine in $\Q$ by \cite{Glazer, Burungale, BhargavaBurungale}. This is the analogue of Theorem \ref{ThmA}.

\item One can ask for the arithmetic analogue of Theorem \ref{ThmB}: is $\Q$ Diophantine in number fields? Our previous work \cite{PastenVidaux} precisely addresses this question and it gives a conditional positive answer.

\item While we were finishing this work, we learned from Arno Fehm that he has arrived at a version of Theorem \ref{ThmA} without the congruence restriction on $q$. His proof is  independent from ours. He obtained further undecidability and definability results for $k(t)$ not covered by the present work. We thank him for sharing his preprint with us.

\end{itemize}
We conclude this introduction by asking whether similar results can be obtained when $k$ is replaced by the algebraic closure of $\F_p$.


\section{The Pheidas predicate} \label{SecPheidas}

This section closely follows our previous work \cite{PastenVidaux}.

Let $K/F$ be a global function field. We introduce the unary predicate $\Theta$ on $K$ defined by
$$
\begin{aligned}
\Theta(b):\quad &\exists \, \nu, \mu, \nu',\mu',y\in K \mbox{ with }y\ne 0\mbox{ and }\\
&\nu^2 + (b-1)\mu^2=b,\quad {\nu'}^2 + (b-1){\mu'}^2=b,\quad \mu\mu'(1-b) = y^2.
\end{aligned}
$$

\begin{lemma} Let $v$ be a place of $K$ and $b\in K$ such that $K\models \Theta(b)$. If $v(b)<0$ then $v(b)$ is even.
\end{lemma}
\begin{proof} This is an elementary verification with valuations that we leave to the reader.
\end{proof}
The following key proposition has a more technical proof, which we leave for Section \ref{SecElliptic}. It is here where we use $q\equiv 3\bmod 4$.
\begin{proposition}\label{PropTheta} Assume $q\equiv 3\bmod 4$. Let $b\in k[t]$ have even degree and $v_t(b)=1$. Then $F\models \Theta(b)$
\end{proposition}
The Pheidas predicate that we use is a unary predicate $\Phi$ on $K$ defined by
$$
\Phi(b): \quad b\in k\mbox{ or }\Theta(t^\epsilon b -a_0 - a_1t) \mbox{ for some }\epsilon \in \{0,1\}\mbox{ and }a_0,a_1\in k.
$$
The subset of $K$ defined by $\Phi$ will also be denoted by $\Phi$.

\begin{theorem}[Properties of the Pheidas predicate] Assume $q\equiv 3\bmod 4$. Then:
\begin{itemize}
\item[(Ph1)] $k[t]\subseteq \Phi$
\item[(Ph2)] There is a finite set of places of $K$, denoted by $S_\Phi$, such that if $b\in \Phi$ then for every $v\notin S_\Phi$ with $v(b)<0$ we have that $v(b)$ is even.
\end{itemize}
\end{theorem}
\begin{proof} This is immediate from the previous results in this section.
\end{proof}


\section{Defining the rational functions} \label{SecField}

Let $n=[K:F]$. For $x\in K$ define its height by $h(x)=n^{-1} \deg \left((x)_\infty\right)$ and
$$
N^{(1)}(x) = \frac{1}{n}\sum_{v_\pfrak(x)>0} \deg \pfrak.
$$ 

\begin{lemma}[The Frobenius relation]\label{LemmaFrob}
The set $\Pcal=\{(u,u^{p^j}):u\in K,\, j\ge0\}$ is Diophantine in $K$. Furthermore, for every $s_0\ge0$ the set $\Fcal=\{(x,x^{q^s},t^{q^s}):x\in K,\, s\ge s_0\}$ is Diophantine in $K$.
\end{lemma}
\begin{proof}
For $\Pcal$ see \cite{PheidasH10ratfn, Shlapentokh}. 

Let $t_0=t^{q^{s_0}}$ and choose $k=\F_p(\zeta)$, with $\zeta=1$ if $k=\F_p$. The conditions $(t_0,T),(t_0+\zeta,T+\zeta)\in\Pcal$ select precisely $T=t^{q^s}$ with $s\ge s_0$. Indeed, comparing degrees in $t_0$ forces the two $p$-powers to agree, and comparing constants forces that power to fix $\zeta$.

To synchronize the power of $x$ and $t$ in order to get $\Fcal$ is routine using $\Pcal$ and the fact that valuation rings are Diophantine \cite{DemeyervanGeel}.
\end{proof}

The next lemma is a standard application of the Riemann--Hurwitz formula for separable maps. For the arithmetic analogue under the $abc$ conjecture, see \cite{PastenVidaux}.
\begin{lemma}[Mason--Stothers]\label{LemmaMason}
Let $g$ be the genus of $K$.  Consider $G\in k[X]$ a squarefree polynomial of degree $d$, and let $x\in K\setminus K^p$. Then 
$$
N^{(1)}(G(x))\ge(d-2)h(x)-\frac{2g-2}{n}.
$$
\end{lemma}
Let $L/K$ be the Galois closure of $K/F$ with Galois group $\Gamma$ and constant field $k_L$. Let us define $f=[k_L:k]$, $m=[L:Fk_L]$, then $|\Gamma| = mf$. Given $s\ge1$ and $\tau\in\Gamma$ acting on $k_L$ by $z\mapsto z^{q^s}$, let $Q_s(\tau)$ be the set of places of $L$ consisting of the unramified primes of $L/F$ with Frobenius $\tau$ above a degree-$s$ prime of $k[t]$.

The following is an application of Chebotarev's theorem in the function field setting, see \cite{FriedJarden}.
\begin{lemma}[Chebotarev]\label{LemmaChebotarev}
There is a constant $c_L$ such that 
$$
\sum_{\Pfrak\in Q_s(\tau)}\deg\Pfrak\ge fq^s-c_Lq^{s/2}
$$ 
for all $s$ and $\tau$ as chosen before.
\end{lemma}


Choose a fixed squarefree $G\in k[X]$ of degree $d\ge 2m+5$ and let $s_0\ge 1$. Define the unary predicate $\Psi(x)$ by requiring that $x\notin K^p$ and that some $w,T,M\in K$ satisfy
\begin{equation}\label{EqnPsi}
(x,w,T)\in\Fcal,\quad (T-t,M)\in\Pcal,\quad \frac{(w-x)^2}{T-t}\in\Phi,\quad \frac{M^2}{G(x)}\in\Phi.
\end{equation}
The predicate $\Psi$ is Diophantine by Lemma \ref{LemmaFrob} and the decomposition $K=\bigoplus_{i=0}^{p-1}K^pt^i$.
\begin{lemma}\label{LemmaPolys}
Every $a\in k[t]- F^p$ satisfies $\Psi$ over $K$.
\end{lemma}
\begin{proof}
Separability of $K/F$ gives $a\notin K^p$. Choose $s\ge s_0$ divisible by the degrees of all the irreducible factors of $G(a)$, and take $p^j\ge\deg G(a)$. For the choice $T=t^{q^s}$, each such factor divides $T-t$, so $G(a)$ divides $(T-t)^{p^j}$. Also $T-t$ divides $a\circ T-a=a^{q^s}-a$. Thus, both fractions in \eqref{EqnPsi} belong to $k[t]\subseteq\Phi$.
\end{proof}
\begin{lemma}\label{LemmaHeight}
There is a constant $c_1$ such that 
$$
(d-4)h(x)\le q^s+c_1
$$ 
whenever $x\in K\setminus K^p$, $s\ge1$, and $j\ge0$ with $M^2/G(x)\in\Phi$, where $M=(t^{q^s}-t)^{p^j}$.
\end{lemma}
\begin{proof}
Away from $S_\Phi$, every odd-order zero of $G(x)$ divides $t^{q^s}-t$: for otherwise $M^2/G(x)$ would have an odd-order pole there. Such primes have total degree $\le nq^s$. Counting each zero of even order with at most half its multiplicity gives
$$
N^{(1)}(G(x))\le\frac{1}{2} h(G(x))+\frac{1}{2} q^s+c_\Phi
=\frac {d}{2} h(x)+\frac{1}{2} q^s+c_\Phi,\quad\mbox{ where } c_\Phi=\frac{1}{n}\sum_{\pfrak\in S_\Phi}\deg\pfrak.
$$
Lemma \ref{LemmaMason} gives the claim with $c_1=2c_\Phi+\frac{4g-4}{n}$.
\end{proof}
%
%
\begin{lemma}\label{LemmaFixed}
If $s_0$ is large enough, every $x\in K$ with $K\models \Psi(x)$ lies in $F$.
\end{lemma}
\begin{proof}
Let $c_S$ be the sum of the degrees of the primes of $L$ above some prime in $S_\Phi$. Since $d-4>2|\Gamma|/f$, we can choose $s_0$ satisfying that for every $s\ge s_0$:
\begin{equation}\label{EqnBound}
2|\Gamma| \frac{q^s+c_1}{d-4}<fq^s-c_Lq^{s/2}-c_S.
\end{equation}
Suppose $\Psi(x)$ holds with $T=t^{q^s}$. Fix $\tau$ acting on $k_L$ by the $q^s$-Frobenius. For $\Pfrak\in Q_s(\tau)$ away from $S_\Phi$, define $\pfrak=\Pfrak\cap K$. Unramifiedness gives $v_\pfrak(T-t)=1$. Hence, the first fraction in \eqref{EqnPsi} has odd valuation, necessarily positive by (Ph2). It follows that $x$ is integral at $\Pfrak$ and $\tau x\equiv x^{q^s}\equiv x\bmod\Pfrak$. If $\tau x\ne x$, the zero divisor of $\tau x-x$ has degree $\ge fq^s-c_Lq^{s/2}-c_S$ by Lemma \ref{LemmaChebotarev}, and its pole divisor has degree at most $2|\Gamma|h(x)$. This contradicts \eqref{EqnBound} and Lemma \ref{LemmaHeight}. The set of elements $\tau$ fixing $x$ (which is a subgroup of $\Gamma$) therefore includes a full coset of $\Gal(L/Fk_L)$ in $\Gamma$, hence, that subgroup too. Thus $x\in Fk_L\cap K=F$, because $K$ has constant field $k$.
\end{proof}
\begin{proof}[Proof of Theorem \ref{ThmB}]
For this $s_0$, Lemmas \ref{LemmaPolys} and \ref{LemmaFixed} give 
$$
k[t] -  F^p\subseteq\Dcal :=\{x:K\models \Psi(x)\}\subseteq F^\times.
$$ 
For nonzero $a,b\in k[t]$, each of the $F^p$-lines $a^{-1}F^p$ and $b^{-1}F^p$ contain at most one of $1,t,t+1$. Let us select $c$ among the three, outside both lines. Then $ca,cb\in\Dcal$, and we deduce that 
$$
F=\{0\}\cup\{u/v:u,v\in\Dcal\}
$$ 
is Diophantine in $K$.
\end{proof}


\section{Defining the polynomials} \label{SecPoly}

In this section $K=F$. For $z\in F^\times$ with $v_\infty(z)\le0$ put $\delta(z)=-v_\infty(z)$; so $\delta(z)=\deg z$ for $z\in A\setminus\{0\}$.
\begin{lemma}[Degrees]\label{LemmaDegrees}
Let $\rho\subseteq\N^r$ be a recursively enumerable relation. Then the set of all $(z_1,\ldots,z_r)\in(F^\times)^r$ with $v_\infty(z_i)\le 0$ for all $i$, and $\rho(\delta(z_1),\ldots,\delta(z_r))$, is Diophantine over $F$.
\end{lemma}
\begin{proof}
Consider the relation $s|_pt$ defined on positive integers by the condition that there is $i\ge 0$ with $t=p^is$. It is a theorem of Pheidas that $P=(\N, 0,1,+,|_p,=)$ positive existentially defines multiplication \cite{PheidasSeries}. Since local rings at a prime are Diophantine in $F$ by \cite{DemeyervanGeel}, using Lemma \ref{LemmaFrob} one gets a positive existential interpretation of $P$ via $\delta$. We conclude by the DPRM property of $\Z$.
\end{proof}

\begin{proof}[Proof of Theorem \ref{ThmA}]
Fix $G$ and $c_1$ as in Section \ref{SecField} with $K=F$. For $s_1,s\ge1$ require $s$ to be divisible by every positive integer up to $q^{s_1}+c_1$. This is a recursively enumerable condition on the pairs $(q^{s_1},q^s)$. Let $\Xi(x)$ be the unary relation on $F$ requiring that $x\notin F^p$, $x$ is integral at the finite primes in $S_\Phi$, and for some $s_1,s$ as before and some $j\ge0$,
\begin{equation}\label{EqnXi}
\frac{(t^{q^{s_1}}-t)^{2p^j}}{G(x)}\in\Phi\mbox{ and } (t^{q^s}-t)x^2\in\Phi.
\end{equation}
This is Diophantine: Lemma \ref{LemmaFrob} defines the required powers, Lemma \ref{LemmaDegrees} imposes the condition on $(\delta(t^{q^{s_1}}),\delta(t^{q^s}))$ and \cite{DemeyervanGeel} gives the necessary integrality.

For $x\in k[t]\setminus F^p$, choose $s_1,j$ as in Lemma \ref{LemmaPolys} and a suitable $s$; both expressions in \eqref{EqnXi} belong to $k[t]$. Conversely, if $\Xi(x)$ holds, Lemma \ref{LemmaHeight} gives $h(x)\le q^{s_1}+c_1$. Any finite pole $\pi$ of $x$ has degree dividing $s$, so, $v_\pi(t^{q^s}-t)=1$. The second expression in \eqref{EqnXi} has a pole of odd order at $\pi$, forcing $\pi\in S_\Phi$, and $x$ is integral. Thus $\Xi$ defines $k[t]\setminus F^p$. We conclude because $k[t]=k\cup \{y^{p^j}:\Xi(y)\mbox{ and }j\ge0\}$ is Diophantine.
\end{proof}


\section{The elliptic curves argument} \label{SecElliptic}

From now on we assume $q\equiv 3 \bmod 4$, so that $-1$ is not a square in $k$. The goal in this section is to prove Proposition \ref{PropTheta}.


\subsection{The two quartics.}

Put $c=1-b$. If $c$ is a square, let us take $\nu=\nu'=\mu=\mu'=1$ and $y^2=c$. Therefore, we can assume from now on that $c$ is not a square.

We claim that for some $\lambda\in F^\times$ the equations
\begin{equation}\label{EqnQuartics}
N^2=bU^4+c\lambda^2,\mbox{ and } {N'}^2=b{U'}^4+c^3\lambda^2
\end{equation}
have solutions with $UU'\ne 0$. Indeed, one sees that the choices
\begin{equation}\label{EqnChoices}
(\nu,\mu)=\left(N/U^2,\lambda/U^2\right),\mbox{ and }(\nu',\mu')=\left(N'/{U'}^2, c\lambda/{U'}^2\right)
\end{equation}
satisfy the definition of $\Theta$, with $y=c\lambda/(UU')$. We prove the claim by constructing $D=b(b-1)\lambda^2$ for which both $D$ and $Dc^2$ satisfy the same quartic-solubility criterion.


\subsection{The rational-point criterion.}
For $D\in F^\times$ we consider the elliptic curves over $F$
$$
E_D:y^2=x^3-Dx\mbox{ and }E'_D:y^2=x^3+4Dx.
$$
Let $\phi:E_D\to E'_D$ and $\hat{\phi}:E'_D\to E_D$ be the isogenies with kernel $\{O,(0,0)\}$.
For $D,d\in F^\times$ let $C_d(D)$ be the quartic 
$$
N^2=dU^4-(D/d)V^4,
$$ 
with $(U,V)\ne(0,0)$. Let us define $S(D)$ as the set of classes $[d]\in F^\times/F^{\times 2}$ for which $C_d(D)$ is everywhere locally soluble, and put $S'(D)=S(-4D)$. From descent theory they are groups.

\begin{lemma}[Rational-point criterion]\label{LemmaRational}
If $D$ and $-D$ are nonsquares and $|S(D)|=4$ and $|S'(D)|=2$, then each $C_d(D)$ with $[d]\in S(D)$ has an $F$-rational point.
\end{lemma}
\begin{proof}
The points $O$ and $(0,0)$ give two classes in each descent image. Thus $S'(D)$ is entirely represented globally, and $S(D)$ has at most one further coset modulo the image of rational points. In terms of the dual 2-isogenies $\phi$ and $\hat{\phi}$ this means $\Sha(E_D)[\phi]=0$ and $|\Sha(E'_D)[\hat{\phi}]|\le 2$. The 2-primary groups $T=\Sha(E_D)[2^\infty]$ and $T'=\Sha(E'_D)[2^\infty]$ are finite by isotriviality \cite{Ulmer}. Hence, $\phi:T\hookrightarrow T'$ because $\Sha(E_D)[\phi]=0$, and $|\ker(\hat{\phi}_*:T'\to T)|\le 2$. Since $T[2]\hookrightarrow\ker\hat{\phi}_*$ has square order by the alternating Cassels--Tate pairing, we get $T[2]=0$, hence $T=0$. Then $T'=\ker\hat{\phi}_*$ also has square order at most $2$, so $T'=0$. The descent exact sequence proves the assertion.
\end{proof}


\subsection{Local conditions.} Let $v$ be a place of $F$. Define $H_v=F_v^\times/F_v^{\times 2}$, $U_v\subset H_v$ the unit classes, and $W_v(D),W'_v(D)$ the local images for $S(D),S'(D)$ respectively. For $D=\varpi^e u$ ($u$ a unit, $\varpi$ a uniformizer at $v$) one computes
\begin{equation}\label{EqnLocTable}
\begin{array}{l | cc}
& W_v(D)&W'_v(D)\\ \hline
e\equiv0\pmod4 & U_v & U_v\\
e\text{ odd} & \langle[-D]\rangle & \langle[D]\rangle\\
e\equiv 2\bmod4,\ u\text{ square} & H_v & \{1\}\\
e\equiv 2\bmod4,\ u\text{ nonsquare} & \{1\}&H_v
\end{array}
\end{equation}
where the last two rows require $-1$ nonsquare in the residue field (because $q\equiv 3\bmod 4$). This follows by comparing $2v(x)$ with $v(D)$ in $y^2=x(x^2-D)$, and using $|W_v(D)|\,|W'_v(D)|=4$ (as $\hat{\phi}\phi = [2]$). For $e\equiv2\bmod 4$ one checks whether $v(x)=1$, with the point $(\varpi\sqrt{u},0)$ available when $u$ is a square. In each case one sees that $W_v(D)$ and $W'_v(D)$ are orthogonal under the Hilbert symbol.


\subsection{The choice of $D$.}
Let $\gamma$ be the monic squarefree representative of $b(b-1)$, and let $B$ be its support. Set $\Omega=\langle[b],[c]\rangle=\langle[b],[-\gamma]\rangle$ and $\Omega'=\langle[\gamma]\rangle$, and note that they have dimensions $2$ and $1$ respectively. Since $\gamma$ is a local square at infinity, we may choose $a\in\{0,1\}$ with $\deg(\gamma)+2a\equiv2\bmod4$ and put $D_0=\gamma t^{2a}$ (here we are using that $b$ has even degree). The table \eqref{EqnLocTable} gives
$$
(W_v(D_0),W'_v(D_0))=
\begin{cases}
(\langle[-\gamma]\rangle,\langle[\gamma]\rangle)& (v\in B),\\
(H_\infty,\{1\})&(v=\infty),\\
(U_v,U_v)&\text{otherwise}.
\end{cases}
$$
Thus $\Omega\subseteq S_0:=S(D_0)$ and $\Omega'\subseteq S'_0:=S'(D_0)$: at a prime factor in $B$ of $b$ (respectively $c$), the element $c$ (respectively $b$) is a local square. These groups are finite; $[-1]$ is excluded by the conditions at $t$ and $\infty$, respectively.

Choose lifts $\delta_1,\ldots,\delta_r$ and $\epsilon_1,\ldots,\epsilon_s$ of bases of $S_0/\Omega$ and $S'_0/\Omega'$. Quadratic reciprocity and Dirichlet's theorem \cite{Rosen} give distinct odd-degree primes $\ell_i$ and $m_j$ outside $B$ with the following quadratic residue properties (here, $[h=i]$ is $1$ if $h=i$ and $0$ otherwise):
\begin{equation}\label{EqnCharacters}
\begin{aligned}
(b/\ell_i)=(c/\ell_i)&=1, &(\delta_h/\ell_i)&=(-1)^{[h=i]},\\
(\gamma/m_j)&=1, &(\epsilon_h/m_j)&=(-1)^{[h=j]}.
\end{aligned}
\end{equation}
Put
\begin{equation}\label{EqnTwist}
D=D_0\left(t^{r+s}\prod_i\ell_i\prod_j m_j\right)^2.
\end{equation}
All the previous local conditions continue to hold, including the prime at infinity, since the factor squared in \eqref{EqnTwist} has even degree. At $\ell_i$, the first image is $\{1\}$ and the second one is $H_{\ell_i}$. On the other hand, at $m_j$ the roles reverse. In particular, $\Omega\subseteq S(D)$ and $\Omega'\subseteq S'(D)$.

On the other hand, if $\eta\in S'(D)$, orthogonality under the Hilbert symbol and Hilbert reciprocity for $(\delta_i,\eta)_v$ imply $1=(\delta_i,\eta)_{\ell_i}=(-1)^{v_{\ell_i}(\eta)}$. Hence, $\eta\in S'_0$, and the conditions at $m_j$ discard all its $\epsilon_j$-coordinates. From a symmetric argument and pairing a class of $S(D)$ with $\epsilon_j$, we deduce
$$
S(D)=\Omega,\qquad S'(D)=\Omega'.
$$
These equalities also hold for $Dc^2$: if $v(c)$ is odd then $v(D)$ is odd, and otherwise multiplication by $c^2$ leaves the corresponding local row of \eqref{EqnLocTable} unchanged. Both parameters $D$ and $Dc^2$ and their negatives are non-squares, since their valuations at $t$ are odd. Writing $D=b(b-1)\lambda^2$, the criterion given by Lemma \ref{LemmaRational} represents $[b]$ on $C_b(D)$ and $C_b(Dc^2)$. Neither point has $U=0$ or $V=0$ (since $[b]\ne[-D],1$), so normalization $V=1$ gives \eqref{EqnQuartics}. The substitutions \eqref{EqnChoices} finish the proof of Proposition \ref{PropTheta}.


\section{Acknowledgements and use of AI}

H.P. was supported by ANID Fondecyt Regular grant 1230507 from Chile. X.V. was supported by Proyecto 2025001477INV from Universidad de Concepci\'on.

The basic strategy of using a Pheidas predicate together with Mason-Stothers comes from our previous work \cite{PastenVidaux} where we used a Pheidas predicate together with the $abc$ conjecture. The idea of using the relation $|_p$ in this context was suggested by AI models. Technical computations were performed with AI assistance and checked by the authors.

We thank Arno Fehm for sharing his related work with us.

\end{document}